\documentclass[12pt,draft]{extarctic}

\usepackage{amsmath,amsthm,amssymb,calc}
\usepackage[dvips]{graphics, graphicx}

\usepackage[english]{babel}

\def\eps{\varepsilon}

\newcounter{num}[section]

\newcommand{\Th}{\refstepcounter{num}
{\bf Theorem \arabic{section}.\arabic{num} }}

\newcommand{\Lemma}{\refstepcounter{num}
{\bf Lemma \arabic{section}.\arabic{num} }}

\newcommand{\Cor}{\refstepcounter{num}
{\bf Corollary \arabic{section}.\arabic{num} }}

\newcommand{\Note}{\refstepcounter{num}
{\it Note \arabic{section}.\arabic{num} }}

\newcommand{\Def}{\refstepcounter{num}
{\it Definition \arabic{section}.\arabic{num} }}

\newcommand{\Proof}{{\bf Proof. }}

\def\eps{\varepsilon}
\def\_phi{\varphi}

\def\a{\alpha}
\def\d{\delta}
\def\la{\lambda}

\def\v{\vec}
\def\F{\widehat}

\def\m{\times}
\def\t{\tilde}
\def\o{\omega}
\def\ov{\overline}

\def\C{{\mathbb C}}

\def\E{{\mathbb E}}

\def\r{\mathcal{R}}
\def\Z_N{{\mathbb Z}_N}
\def\Z{{\mathbb Z}}
\def\N{{\mathbb N}}
\def\U{{\mathcal U}}

\def\Gr{{\mathbf G}}

\def\cs{{\rm cs\,}}

\def\D{{\mathbb D}}

\def\l{\left}
\def\r{\right}

\author{Shkredov I.D.}
\title{ On Gowers norms of some functions
\footnote{ This work was supported
Pierre Deligne's grant based on his 2004 Balzan prize, grant RFFI
N 06-01-00383, grant of President of Russian Federation
MK-1959.2009.1, and grant Leading Scientific Schools No.
691.2008.1.
\newline
{\bf Keywords} : Gowers norms, linear equations.
\newline
MSC 2000 : 11B75, 11B99.} }
\date{}
\begin{document}
\maketitle

\begin{center}
 Annotation.
\end{center}

{\it \small
    We consider a class of two--dimensional functions $f(x,y)$ with the property that
    the smallness of its rectangular norm implies the smallness of
    rectangular norm for $f(x,x+y)$.
    Also we study a family of functions $f(x,y)$
    having a similar property
    for higher Gowers norms.
    The method based on a transference principle for a class of sums
    over special systems of linear equations.
}
\\
\\
\\

\refstepcounter{section} \label{sec:introduction}

\begin{center}
{\bf \arabic{section}. Introduction.}
\end{center}

    The notion of Gowers norms was introduced in papers \cite{Gow_4,Gow_m}
    and is a very important tool of investigation in wide class
    of problems of additive combinatorics
    (see e.g. \cite{Gow_4,Gow_m} \cite{GT_great}---\cite{Shkr_abealian})
    as well as in ergodic theory (see e.g. \cite{TZ}, \cite{Fu}---\cite{Austin_norm}).
    Recall the definitions.

  Let
  $\Gr$
  be a finite
  set,
  and $N=|\Gr|$.
Let also $d$ be a positive integer, and $\{ 0,1 \}^d = \{ \omega =
(\omega_1,\dots, \omega_d) ~:~ \omega_j \in \{0,1\}, j=1,2,\dots,d
\}$ be an ordinary $d$---dimensional cube. For $\omega \in
\{0,1\}^d$ denote by $|\omega|$ the sum $\omega_1 + \dots +
\omega_d$.
Let also $\mathcal{C}$ be the operator of complex conjugation. Let
$\v{x} = (x_1,\dots,x_d), \v{x}'=(x'_1,\dots,x'_d)$ be two
arbitrary vectors from $\Gr^d$. By $\v{x}^\o = (\v{x}^\o_1, \dots,
\v{x}^\o_d)$ denote the vector
\begin{displaymath}
    \v{x}^\o_i = \left\{ \begin{array}{ll}
                        x_i & \mbox{ if } \o_i = 0, \\
                        x'_i & \mbox{ if } \o_i = 1. \\
                        \end{array} \right.
\end{displaymath}
Thus $\v{x}^\o$ depends on $\v{x}$ and $\v{x}'$.

Let $X$ be a non--empty finite set, $Z : X \to \C$ be a function.
Denote by $\mathbb{E} Z = \mathbb {E}_x Z$ the sum $\frac{1}{|X|}
\sum_{x\in X} Z(x)$.

    Let $f : \Gr^d \to \C$ be an arbitrary function.
    We will write $f(\v{x})$ for $f(x_1,\dots,x_d)$.

\Def {\it (see \cite{Gow_4,Gow_m})} {
    {\it Gowers $U^d$--norm} (or $d$--uniformity norm) of the function $f$ is the following expression
\begin{equation}\label{f:G_norm_d'}
    \| f \|_{U^d}
            =
                \left( \E_{\v{x}\in \Gr^d}\, \E_{\v{x}' \in \Gr^d}
                    \prod_{\o \in \{ 0,1 \}^d} \mathcal{C}^{|\o|} f (\v{x}^{\o}) \right)^{1/2^d} \,.
\end{equation}
}

A sequence of $2^d$ points $\v{x}^\o$, $\o \in \{ 0,1 \}^d$ is
called {\it $d$--dimensional cube}. Thus the summation in formula
(\ref{f:G_norm_d'}) is taken over all cubes of $\Gr^d$. For
example, $\{ (x,y), (x',y), (x,y'), (x',y') \}$, where
$x,x',y,y'\in \Gr$ is a two--dimensional cube. In the case Gowers
norm is called {\it rectangular norm}.

For $d=1$ the expression above gives a semi--norm but for $d\ge 2$
Gowers norm is a norm. In particular, the triangle inequality
holds
\begin{equation}\label{e:triangle}
    \| f+g \|_{U^d} \le \| f \|_{U^d} + \| g \|_{U^d} \,.
\end{equation}
One can prove also (see \cite{Gow_m}) the following monotonicity
relation. Let $f_{x_d} (x_1, \dots, x_{d-1}) := f (x_1, \dots,
x_{d})$. Then
\begin{equation}\label{e:Gowers_mon}
    \E_{x_d \in \Gr} \| f_{x_d} \|^{2^{d-1}}_{U^{d-1}} \le \| f \|^{2^{d-1}}_{U^d}
\end{equation}
for all $d\ge 2$.

If $\Gr = (\Gr,+)$ is a finite Abelian group with additive group
operation $+$, $N=|\Gr|$ then one can "project"\, the norm above
onto the group $\Gr$ and obtain the ordinary Gowers norm. In other
words, we put the function $f(x_1,\dots,x_d)$ in formula
(\ref{f:G_norm_d'}) equals "one--dimensional"\, function
$f(x_1,\dots,x_d) := f(x_1+\dots+x_d)$.
Denoting the obtained norm as $U^k$ and we have an analog of
(\ref{e:Gowers_mon})
\begin{equation}\label{e:Gowers_mon_1}
    \| f \|_{U^{d-1}} \le \| f \|_{U^d}
\end{equation}
for all $d\ge 2$.

Gowers norms have the following characteristic property. Let $d\ge
2$, $\v{x} = (x_1, \dots, x_d) \in \Gr^d$ be an arbitrary vector
and $i\in \{1,\dots,d \}$. By $(\v{x})_{(i)}$ denote the vector
$(\v{x})_{(i)} = (x_1, \dots, x_{i-1}, x_{i+1}, \dots, x_d) \in
\Gr^{d-1}$. Applying the Cauchy--Schwartz (see \cite{Gow_m} and
\cite{Gow_hyp4,Gow_hypk}) several times, we have

\Lemma {\it
    Let $d\ge 2$ be an integer, and $f:\Gr^k \to \C$ be a function.
    Let also $u_1, \dots, u_d : \Gr^d \to [-1,1]$ be any functions such that
    $u_i (\v{x}) = u_i ((\v{x})_{(i)})$, $i=1,\dots,d$, $\v{x} = (x_1,\dots,x_d)$.
    Then
    $$
        \l| \sum_{\v{x}} f(\v{x}) \prod_{i=1}^d u_i (\v{x}) \r|
                \le
                    \| f \|_{U^d} \,.
    $$
} \label{l:U_k_sense}

Thus any function with small Gowers norm does not correlate with
product of any functions which depend on smaller number of
variables.

In the paper we
concentrate on the case of two--dimensional functions $f:\Gr \m
\Gr \to \C$, $\Gr$ is a finite Abealian group. For any positive
integer $t$
one can consider $U^{t+1}$---norm of the function $F_{t}
(x,y_1,\dots,y_t) : = f(x, x+y_1+\dots+y_t)$ and the function
$H_{t} (x,y_1,\dots,y_t) : = f(x, y_1+\dots+y_t)$.
It is easy to construct examples of functions $f(x,y)$ with, say,
huge rectangular norm and small quantity $\| f(x,x+y) \|_{U^2}$
("skew rectangular norm") and vice versa. Thus there is no obvious
dependence between the numbers $\| F_{t} \|_{U^{t+1}}$ and $\|
H_{t} \|_{U^{t+1}}$ in general.
Nevertheless, we find a class of functions such that the smallness
of $\| H_{t} \|_{U^{t+1}}$ implies the smallness $\| F_t
\|_{U^{t+1}}$.

By $\F{\Gr}$ denote the Pontryagin dual of $\Gr$. In other words
$\F{\Gr}$ is the group of homomorphisms $\xi$ from $\Gr$ to
$\mathbf{R} / \mathbf{Z}$, $\xi : x \to \xi \cdot x$. It is
well--known that in the case of Abelian group $\Gr$ the dual group
$\F{\Gr}$ is isomorphic to $\Gr$.

Let us formulate one of the main results of the paper.

\Th {\it
    Let $a:\Gr \to \Gr$ be a function and
    $f (x,y) = e(x \cdot a(y))$,
    where $e(x) = e^{2\pi i x}$.
    Then the condition
    \begin{equation}\label{f:intr_1}
        \| f(x,y_1+\dots+y_t) \|_{U^{t+1}} = o (1) \,, \quad N\to \infty
    \end{equation}
    implies
    \begin{equation}\label{f:intr_2}
        \| f(x, x+y_1+\dots+y_t) \|_{U^{t+1}} = o (1) \,, \quad N\to \infty
    \end{equation}
    for all positive integers $t$.

} \label{t:intr_rect}

More precisely, we obtain a quantitative form of the formulas
above (see Corollary \ref{cor:rectangle_imply}). Note that more
strong condition than (\ref{f:intr_1}), namely, $\|
f(x,y_1+\dots+y_{t+1}) \|_{U^{t+1}} = o(1)$ trivially  implies
(\ref{f:intr_2}) (see proof of Corollary
\ref{cor:rectangle_imply}). Theorem \ref{t:intr_rect} will be
derived from the more general Theorem \ref{t:main_skew_norms} of
section \ref{sec:proof}. The result
says, roughly speaking, that
a class of sums of general form is taken over systems of linear
equations including an arbitrary function $a(x)$ can be reduced to
sums having no these linear restrictions (see Theorem
\ref{t:main_skew_norms}). It is interesting we do not use Fourier
analysis in the proof. Our main tool is a suitable version of
Lemma 9.3 from Gowers' paper \cite{Gow_m}.

Let us say a few words about the notation. If $S\subseteq \Gr$ is
a set then we will write $S(x)$ for the characteristic function.
In other words $S(x) = 1$ if $x\in S$ and zero otherwise. By
$\log$ denote logarithm base two and by $\D$ denote the unit disk
on the complex plane. Sings $\ll$ and $\gg$ are usual Vinogradov's
symbols. If $n$  is a positive integer then we will write $[n]$
for the segment $\{ 1,2,\dots,n \}$.

The author is grateful to N.G. Moshchevitin for useful
discussions.

\refstepcounter{section}
\label{sec:proof}

\begin{center}
{\bf \arabic{section}. The proof of the main result.}
\end{center}

We need in a number of definitions.

Let $l\ge 2$, $m$ be any positive integers, $m<l$.
Let also $x_1,\dots,x_l$, $y_1,\dots,y_l$ be some variables.
Consider a system of linear equations
\begin{equation}\label{eq:1}
    \sum_{j=1}^l \eps^{(i)}_j x_j = 0 \,,\quad i=1,2,\dots, m \,,
\end{equation}
and also an equation
\begin{equation}\label{eq:2}
    \sum_{j=1}^l \eps^{(1)}_j y_j = 0 \,.
\end{equation}
Here $\eps^{(i)}_j \in \{ 0,-1,1 \}$.
So, if we know system  (\ref{eq:1}) then we automatically
know equation (\ref{eq:2}).
Suppose that the rank of subsystem (\ref{eq:1}) equals $m$.

Let $S$ be the  family of all systems from $m$ equations of full rank
having the form
$\sum_{j=1}^l \eta_j x_j = 0$, $\eta_j \in \{ 0,\pm 1\}$ such that $\eta$ can
be written as a linear combination of vectors
$(\eps^{(1)}_1,\dots, \eps^{(1)}_l),\dots,(\eps^{(m)}_1,\dots, \eps^{(m)}_l)$
from (\ref{eq:1}).
For example, we can multiply some equations of (\ref{eq:1}) by $-1$
and get a system from $S$.
Thus $|S| \ge 2^m$.
On the other hand, $|S| \le 3^{lm}$.
Clearly, if a tuple $(x_1,\dots,x_l)$ is a solution of some system from $S$
then it is a solution of any system from $S$, in particular, system (\ref{eq:1}).
By $\Upsilon$ denote system (\ref{eq:1}), (\ref{eq:2}).

Let $E$ be a family of linear equations with coefficient $\{0,\pm 1\}$
which can be written as a linear combinations of some equations from $S$.
Clearly, $|E| \le 3^l$.
Suppose that $e\in E$ and define $\theta(e)=2^{t}$,
where $t$ is the number of zero coefficients in $e$.
Put also $\theta(\upsilon) = \prod_{e\in \upsilon} \theta(e)$,
$\upsilon$ is a system.
Lastly, let
$$
    \theta = \theta_{m-1} (E) = \sum_{\upsilon} \theta(\upsilon) \,,
$$
where summation is taken over all systems with $m-1$ equations from $E$.

Let $B_1,\dots,B_l \subseteq \Gr$ be
arbitrary
sets.
Suppose that there are some maps $\_phi_1,\dots,\_phi_l$, $\_phi_j : B_j \to \Gr$, $j\in [l]$.

Finally, let $C$ be a subset of $\Gr^l$.
We will call a tuple $(x_1,\dots,x_l)$ satisfying
a system from $S$
{\it additive}.
An additive tuple $(x_1,\dots,x_l)$ is called {\it degenerate} if there is a non--zero vector
$\v{\eta} = (\eta_1,\dots, \eta_l)$, $\eta_j \in \{ 0,\pm 1\}$
such that $\sum_{j=1}^l \eta_j x_j = 0$ and the last equation does not
belong $E$.
Otherwise, such a tuple will be called {\it non-degenerate}.
Note that non--degenerate additive tuples can satisfy equations from $E$ only.
Further, non-degenerate additive tuple $(x_1,\dots,x_l) \in C$
is called {\it good} if
$(x_1,\dots,x_l)$, $(\_phi_1 (x_1),\dots,\_phi_l (x_l))$ satisfies
(\ref{eq:1}), (\ref{eq:2})
for some system from $S$,
and {\it bad} otherwise.
So, by definition, good and bad tuples are non--degenerate and belong to $C$.

First of all we prove the following simple extension of Lemma 9.3 from \cite{Gow_m}.

\Lemma
{\it
    Let $\a,\o,\eta \in (0,1]$ be any numbers,
    $B_j \subseteq \Gr$ be arbitrary sets,
    and $C$ be a subset of $\Gr^l$.
    Suppose that  $T$ is a parameter,
    the number of additive tuples in $(B_1 \m \dots \m B_l) \bigcap C$
    is
    at most
    $\o^{-1} T$,
    \begin{equation}\label{cond:N}
         T \ge 2^{2l} |S| \l( \frac{2(1+\o^{-1})}{\a \eta} \r)^{lm 2^{3ml}} N^{l-m-1}
            \,,
    \end{equation}
    and
    there exist
    at least
    $\a T$
    good tuples.
    Then there are some sets $B'_j \subseteq B_j$
    such that the number of good tuples
    $(x_1,\dots,x_l) \in (B'_1 \m \dots \m B'_l) \bigcap C$
    is at least
    $$
        \l( \frac{\a \eta}{2(1+\o^{-1})} \r)^{lm 2^{3ml}} T 
    $$
    and
    the ratio of good tuples
    in $(B'_1 \m \dots \m B'_l) \bigcap C$
    to the number of bad tuples
    in $(B'_1 \m \dots \m B'_l) \bigcap C$
    is at least $(1-\eta)$.
}
\label{l:Gowers_ext}
\\
\Proof
Let $k$ be a natural parameter and
choose $r_1,\dots,r_k$, $s_1,\dots,s_k$,
$w^{(1)}_1,\dots,w^{(1)}_k, \dots, w^{(m-1)}_1,\dots,w^{(m-1)}_k \in \Gr$
uniformly and independently.
After that independently choose points $x_j \in B_j$ such that
a point $x_j$ go into a new set $B'_j$ with probability
$$
    p(x_j)
        =
            \frac{1}{2^{mk}} \prod_{i=1}^k \l( 1+\cs (r_i x_j + s_i \_phi_j (x_j)) \r)
                \cdot
            \prod_{q=1}^{m-1} \l( 1+\cs (w^{(q)}_i x_j) \r) \,,
$$
where $\cs (x) = \cos (\frac{2\pi}{N} x)$.
A tuple $(x_1,\dots,x_l)$
belongs to $B'_1 \m \dots \m B'_l$ with probability
$$
    \frac{1}{2^{mlk} N^{(m+1)k}} \sum_{r_1,\dots,r_k}\, \sum_{s_1,\dots,s_k}\,
        \sum_{w^{(1)}_1,\dots,w^{(1)}_k,\dots, w^{(m-1)}_1,\dots,w^{(m-1)}_k}\,
            \prod_{i=1}^k \prod_{j=1}^l
                \l( 1+\cs (r_i x_j + s_i \_phi_j (x_j)) \r)
                    \m
$$
$$
                    \m
                    \prod_{q=1}^{m-1} \l( 1+\cs (w^{(q)}_i x_j) \r)
                        =
$$
$$
                        =
    \frac{1}{2^{mlk} N^{(m+1)k}}
        \l( \sum_{r,s} \sum_{w_1,\dots,w_{m-1}} \prod_{j=1}^l
            (1+\cs (r x_j + s \_phi_j (x_j))) \prod_{q=1}^{m-1} ( 1+\cs (w_q x_j) ) \r)^k
                =
$$
$$
    = \frac{1}{2^{mlk} N^{(m+1)k}}
        ( 2^{-lm} \sum_{r,s} \sum_{w_1,\dots,w_{m-1}}
$$
$$
                \prod_{j=1}^l
            (2+e(r x_j + s \_phi_j (x_j)) + e(-r x_j - s \_phi_j (x_j)))
                \prod_{q=1}^{m-1} ( 2 + e(w_q x_j) + e(-w_q x_j) ) )^k
                    =
$$
$$
    = \frac{1}{2^{mlk} N^{(m+1)k}}
        ( 2^{-lm} \sum_{r,s} \sum_{w_1,\dots,w_{m-1}} \prod_{j=1}^l
        \sum_{\eps^{(1)}_1, \dots, \eps^{(1)}_l,\dots, \eps^{(m)}_1, \dots, \eps^{(m)}_l \in \{0,\pm 1\}}
            2^{|\{ i,j ~:~ \eps^{(i)}_j = 0 \}|}
                    \m
$$
$$
                    \m
                e(r \sum_j \eps^{(1)}_j x_j + s \sum_j \eps^{(1)}_j \_phi_j (x_j) +
                    w_1 \sum_j \eps^{(2)}_j x_j + \dots + w_{m-1} \sum_j \eps^{(m)}_j x_j ) )^k \,.
$$
Thus, the last probability does not equal zero if there is
a tuple $(\eps^{(i)}_j)$, $\eps^{(i)}_j \in \{ 0,\pm 1\}$ such that
\begin{equation}\label{f:eps_ij}
    \sum_j \eps^{(1)}_j x_j
        =
            \sum_j \eps^{(1)}_j \_phi_j (x_j)
                =
                    \sum_j \eps^{(2)}_j x_j
                        =
                            \dots
                                =
                                    \sum_j \eps^{(m)}_j x_j
                                        =
                                            0 \,.
\end{equation}
It is so, for example, if all $\eps^{(i)}_j$ equal zero.

%
%
%
    Let $(x_1,\dots,x_l)$ be a bad tuple.
    Then, by definition, $(x_1,\dots,x_l)$ satisfy a system $\upsilon \in S$
    and is non--degenerate.
    It is easy to see that there is no vector $(\eps^{(1)}_1, \dots, \eps^{(1)}_l) \neq \v{0}$
    such that (\ref{f:eps_ij}) holds.
    Indeed, otherwise we can add some equations from $\upsilon$
    and obtain a contradiction with the fact that $(x_1,\dots,x_l)$ is a bad sequence.
    Thus, every bad additive tuple is chosen with probability
    $2^{-lmk} (2^{-lm} 2^l \theta)^k$.
    Clearly, $\theta \ge 2^{l(m-1)}$ and put $\theta = 2^{l(m-1)} + \theta_1$, $\theta_1 \ge 0$.
    Then the last probability equals
    $2^{-lmk} (1 + \theta_1 2^{-l(m-1)})^k$.

On the other hand, it is easy to see that
a good additive tuple is chosen with probability
at least
$$
    2^{-lmk} (2^{-lm} (2^{l} \theta + 2 \theta) )^k
        \ge 2^{-lmk} (1 + \theta_1 2^{-l(m-1)} + 2^{-(l-1)})^k
            > 2^{-lmk} (1 + \theta_1 2^{-l(m-1)})^k \,.
$$
Note that there are at most
$
 3^l |S| |B_1| \dots |B_{l-m-1}|
    \le
        3^{l} |S| N^{l-m-1}
$
degenerate tuples.
Now, let $X$ and $Y$ be the numbers of good and bad additive tuples in
$B'_1 \m \dots \m B'_l$.
Using the assumption of the lemma, we get
$$
    \E X \ge 2^{-lmk} (1 + \theta_1 2^{-l(m-1)} + 2^{-(l-1)})^k \a T
$$
and
$$
    \E Y \le 2^{-lmk} (1 + \theta_1 2^{-l(m-1)})^k \o^{-1} T \,.
$$
It is easy to see that $\theta_1 < \theta \le 3^{lm} 2^{lm}$.
Since
$(1+2^{-3lm})^k \ge 2^{k 2^{-3lm}}$ it follows that
$$
    \l( \frac{1 + \theta_1 2^{-l(m-1)} + 2^{-(l-1)}}{1 + \theta_1 2^{-l(m-1)}} \r)^k
        >
            (1+2^{-3lm})^k
                \ge
                    2^{k 2^{-3lm}}
$$
and the last
expression is at least $(1+\o^{-1})/(\a \eta)$, provided by we choose an integer $k$ from the conditions
\begin{equation}\label{tmp:3}
    2\l( \frac{\a\eta}{2(1+\o^{-1})} \r)^{2^{3lm}} \le 2^{-k} \le \l( \frac{\a\eta}{1+\o^{-1}} \r)^{2^{3lm}} \,.
\end{equation}
Hence
\begin{equation}\label{tmp:4}
    \eta \E X - \E Y \ge \a \eta 2^{-lmk} (1 + \theta_1 2^{-l(m-1)} + 2^{-(l-1)})^k T
        -
            2^{-lmk} (1 + \theta_1 2^{-l(m-1)})^k \o^{-1} T
                \ge
                    2^{-lmk} T
                    \,.
\end{equation}
Using (\ref{cond:N}), (\ref{tmp:3}),
we have
\begin{equation}\label{tmp:5}
    2^{-lmk} T
        \ge
            2 \l( \frac{\a\eta}{2(1+\o^{-1})} \r)^{lm 2^{3ml}} T
                \ge
                    2 \cdot 4^l |S| N^{l-m-1}
                        \ge
                            2 (3^l |S| |B_1| \dots |B_{l-m-1}|) \,.
\end{equation}
By (\ref{tmp:4}), (\ref{tmp:5}), we obtain that there are sets
$B'_1 \subseteq B_1, \dots, B'_l \subseteq B_l$ such that
$\eta X \ge Y$ and
$X\ge \l( \frac{\a\eta}{2(1+\o^{-1})} \r)^{lm 2^{3ml}} T$.
This completes the proof.

\Note
Certainly, one can suppose that our tuples
$(x_1,\dots,x_l)$, $(\_phi(x_1),\dots, \_phi(x_l))$, $x_j\in B'_j$
satisfying (\ref{eq:1}) and (\ref{eq:2}) for some system $\upsilon$ from $S$,
satisfy a new system
\begin{equation}\label{eq:1'}
    \sum_{j=1}^l \bar{\eps}^{(i)}_j x_j = 0 \,,\quad i=1,2,\dots, m \,,
\end{equation}
and
\begin{equation}\label{eq:2'}
    \sum_{j=1}^l \bar{\eps}^{(1)}_j \_phi(x_j) = 0 \,,
\end{equation}
where $\bar{\eps}^{(i)}_j \in \{0,\pm 1\}$,
$\bar{\eps}^{(i)}_j = \eps^{(i)}_j$ for any $j\ge 1$, $i\ge 2$,
the rank of system (\ref{eq:1'}) equals $m$,
and $\bar{\eps}^{(1)}_j$ depend on $\upsilon$.

\Def
{
Let $L,l$ be positive integers.
A set $\Omega \subseteq \Gr^l$ is called a {\it set of level $L$} if
$$
    \Omega = \prod_{j=1}^L (\Omega'_j-\Omega''_j)
$$
where each $\Omega'_j$, $\Omega''_j$, $\Omega''_j \subseteq \Omega'_j$ is a
Cartesian product of some sets.
We
say that Cartesian products are sets of level zero.
}


Let $C$ be a set, and let $f$ be a complex function which depends on variables
$(x_1,\dots,x_l) \in (B_1\m \dots \m B_l) \cap C$, $\v{x} = (x_1,\dots,x_l)$
and $f$ also depends on $\_phi_1 (x_1),\dots, \_phi_l (x_l)$,
$\_phi (\v{x}) = (\_phi_1 (x_1),\dots, \_phi_l (x_l))$.
Let $\upsilon$ be a system from $S$ and
\begin{equation}\label{f:sigma_form}
    \sigma_{f,\upsilon} (B_1,\dots,B_l ; C) :=
        \sum_{(\v{x},\_phi(\v{x})) \mbox{ satisfies } (\ref{eq:1}), (\ref{eq:2})
            \mbox{ for } \upsilon \in S }
            f(x_1,\dots,x_l,\_phi_1 (x_1),\dots, \_phi_l (x_l))
            \,.
\end{equation}
Further
$
 \sigma_{f} (B_1,\dots,B_l ; C)
    =
        \sum_{\upsilon \in S} |\sigma_{f,\upsilon} (B_1,\dots,B_l ; C)|
$.
Take an arbitrary number $\rho \in [l]$,
express $x_\rho$ from (\ref{eq:1}), $\_phi_\rho (x_\rho)$ from (\ref{eq:2})
and substitute $x_\rho$, $\_phi_\rho (x_\rho)$ into $f$.
After that choose another $m-1$ variables $x_{j_1},\dots,x_{j_{m-1}}$
express them from (\ref{eq:1}) and substitute into $f$.
Let $\v{j} = (j_1,\dots,j_{m-1})$.
We get a new function
$$
    f^\upsilon_{\rho,\v{j}} (\v{x}) :=
    f^\upsilon_{\rho,\v{j}}
    (x_j, \_phi_j (x_j))\,, \quad \rho \notin \{ j_1,\dots,j_{m-1} \} \,.
$$
Let $\| f \|_{\U^d} = N^{2d} \| f \|^{2^d}_{U^d}$ and for any
$\rho_1,\rho_2 \in [d]$, $\rho_1<\rho_2$, we put
$$
    \| f \|_{\U^2 (\rho_1,\rho_2)}
        =
            \sum_{x_1,\dots,x_d}\, \sum_{x'_{\rho_1}, x'_{\rho_2}}
                f(x_1,\dots,x_{\rho_1},\dots,x_{\rho_2},\dots,x_d)
                \ov{f(x_1,\dots,x'_{\rho_1},\dots,x_{\rho_2},\dots,x_d)}
                    \m
$$
$$
                    \m
                \ov{f(x_1,\dots,x_{\rho_1},\dots,x'_{\rho_2},\dots,x_d)}
                f(x_1,\dots,x'_{\rho_1},\dots,x'_{\rho_2},\dots,x_d) \,.
$$

Formulate the main result of the section.

\Th
{\it
    Let $\eps \in (0,1]$ be a real number,
    $\eps \le 2^{-2^{20} l^2 m^2 2^{6ml}}$,
    and $f$ be a complex function, $f:\Gr \to \D$.
    Let also $B_1,\dots,B_l$ be arbitrary sets.
    Suppose that for any $\upsilon \in S$
    there are $\v{j}$,
    $\rho \in [l]$ and
    $\rho_1\neq \rho,\rho_2 \neq \rho$,
    $\rho_1 \neq \rho_2$
    such that
    \begin{equation}\label{cond:rectangle_small}
        \| f^\upsilon_{\rho,\v{j}} (\v{x}) \|_{\U^2 (\rho_1,\rho_2)}
            \le
                \eps N^{l-m+2} \,.
    \end{equation}
    Then
    \begin{equation}\label{concl:sigma_small}
        \sigma_f (B_1,\dots,B_l ; \Gr^l)
            \le
                \max\l\{ 2^{30} \l( \frac{512}{\log (1/\eps)} \r)^{(128 lm 2^{3ml})^{-1}},
                        4 \l( \frac{2^{2l}}{N} \r)^{(16lm 2^{3ml})^{-1}} \r\}
                3^{lm} N^{l-m} \,.
    \end{equation}
}
\label{t:main_skew_norms}
\\
\Proof
Let $\eps_1$ be the
maximum
in the right hand side of (\ref{concl:sigma_small}) divided by $N^{l-m}$.
Let also $\upsilon\in S$ and denote by $\tau_{f,\upsilon} (B_1,\dots,B_l ; C)$ the number
of solutions of system (\ref{eq:1'}), (\ref{eq:2'}).
By $\tau_{f} (B_1,\dots,B_l ; C)$ denote the sum
$\sum_{\upsilon \in S} \tau_{f,\upsilon} (B_1,\dots,B_l ; C)$.
Suppose that inequality (\ref{concl:sigma_small}) does not hold.
Since $\tau_{f} (B_1,\dots,B_l ; \Gr^l) \le |S| N^{l-m}$ and
\begin{equation}\label{f:sigma_le_tau}
    \sigma_{f} (B_1,\dots,B_l ; \Gr^l) \le \tau_{f} (B_1,\dots,B_l ; \Gr^l) \,,
\end{equation}
it follows that
$\sigma_{f} (B_1,\dots,B_l ; \Gr^l) \ge \eps_1 |S|^{-1} \tau_{f} (B_1,\dots,B_l ; \Gr^l)$
and
$\tau_{f} (B_1,\dots,B_l ; \Gr^l) \ge \eps_1 N^{l-m}$.

We need in a "density increment"\, lemma.

\Lemma
{\it
        Let $\a,\o \in (0,1]$ be any numbers,
    $B_j$ be arbitrary sets,
    and
    $C$ be a set of level $L$.
    Suppose that
    \begin{equation}\label{cond:N_new}
         \tau_{f} (B_1,\dots,B_l; C)
            \ge
                |S| \max\{ 2^{2l} \l( \frac{32(1+\o^{-1})|S|}{\a \eps_1} \r)^{lm 2^{3ml}} N^{l-m-1}, \o N^{l-m} \}\,,
    \end{equation}
    \begin{equation}\label{cond:eps_small}
        \eps
            \le
                2^{-4L-16} |S|^{-8} \eps^4_1 \l( \frac{\a \eps_1}{32(1+\o^{-1})|S|} \r)^{4lm 2^{3ml}}
                    \l( \frac{\tau_{f} (B_1,\dots,B_l; C)}{N^{l-m}} \r)^4
    \end{equation}
    and $f : \Gr \to \D$ is a
    function  satisfying (\ref{cond:rectangle_small}).
    Let also
    \begin{equation}\label{f:sigma_a}
        \sigma_{f} (B_1,\dots,B_l; C) \ge \a \tau_{f} (B_1,\dots,B_l; C) \,.
    \end{equation}
    Then
    there is
    a set $\t{C} \subseteq C$ of level at most $L+1$
    and there are sets $B'_j$, $j\in [l]$,
    $(B'_1 \m \dots \m B'_l) \bigcap \t{C} = \emptyset$
    such that
    \begin{equation}\label{f:sigma_a+}
        \sigma_{f} (B_1,\dots,B_l; C)
            \le
                \sigma_{f} (B_1,\dots,B_l; \t{C})
                            +
                                2^{-2} \eps_1 |S|^{-1} \tau_f (B'_1,\dots,B'_l ; C) \,.
    \end{equation}
    and
    \begin{equation}\label{f:solutions_DONOT_decrease}
        \tau_{f} (B_1,\dots,B_l; \t{C})
            \le
                (1-\zeta) \tau_{f} (B_1,\dots,B_l; C) \,,
    \end{equation}
    where
    $\zeta = \l( \frac{\a \eps_1}{32(1+\o^{-1})|S|} \r)^{lm  2^{3ml}}$.
}
\label{l:density_inc}
\\
\Proof
Using Lemma \ref{l:Gowers_ext} with parameters
$\a,\o,\eta = |S|^{-1} \eps_1/16$
and $T=\tau_{f} (B_1,\dots,B_l; C)$,
we get the sets $B'_j \subseteq B_j$, $j\in [l]$
such that
the number $g$ of good tuples
    $(x_1,\dots,x_l) \in (B'_1 \m \dots \m B'_l) \bigcap C$
    is at least
    \begin{equation}\label{f:number_of_good_tuples_below}
        g \ge \l( \frac{\a \eta}{2(1+\o^{-1})} \r)^{lm 2^{3ml}}
        \tau_{f} (B_1,\dots,B_l; C)
            =
                \zeta \tau_{f} (B_1,\dots,B_l; C) \,.
    \end{equation}
Thus
\begin{equation}\label{tmp:10.12.2009_1}
    \tau_f (B'_1,\dots,B'_l ; C) \ge g \ge \zeta \tau_f (B_1,\dots,B_l ; C) \,.
\end{equation}
Let $(B'_j)^1 = B'_j$ and $(B'_j)^0 = B_j \setminus B'_j$.
Then for any $\upsilon \in S$, we have
$$
    \sigma_{f,\upsilon} (B_1,\dots,B_l ;C)
        = \sum_{\varpi \in \{0,1\}^l} \sigma_{f,\upsilon} ( (B'_1)^{\varpi_1}, \dots, (B'_l)^{\varpi_l};C)
            = \sum_{\varpi \in \{0,1\}^l} \sigma_{f,\upsilon} (\varpi;C) \,,
$$
where $\varpi = (\varpi_1,\dots, \varpi_l)$.
Without losing of generality, assume that
$\rho=l$ and $\v{j}=(x_{l-m+1}, \dots, x_{l-1})$.
Consider the term
$\sigma_1 (\upsilon):= \sigma_{f,\upsilon} ( B'_1, \dots, B'_l; C)$
which corresponds to $\varpi_1 = \dots = \varpi_l = 1$.
We have
$$
    \sum_{\upsilon \in S} |\sigma_1 (\upsilon)|
        \le
              \sum_{\upsilon \in S}
              \l| \sum_{\v{x} \in (B'_1\m \dots \m B'_{l}) \bigcap C,\,
                    (\v{x},\_phi(\v{x})) \mbox{ satisfies }
                        (\ref{eq:1}), (\ref{eq:2}) \mbox{ for } \upsilon}
                      f^\upsilon_{\rho,\v{j}} (\v{x}) \r|
                    +
$$
\begin{equation}\label{21.10.2009_1}
                    +
                        \eta \tau_f (B'_1,\dots,B'_l ; C)
                            =
                                \sigma'_1 + \sigma''_1
                                \,.
\end{equation}
It is easy to see, using
(\ref{cond:rectangle_small}) and Lemma \ref{l:U_k_sense}
that
\begin{equation}\label{21.10.2009_2}
    \sigma'_1 \le 2^L |S| \eps^{1/4} N^{l-m}
        \le
            2^L 3^{lm} \eps^{1/4} N^{l-m} \,.
\end{equation}

Put
$\t{C} = ((B_1\m \dots \m B_l) \setminus (B'_1\m \dots \m B'_l)) \bigcap C$.
Clearly, $\t{C}$ is a set of level $L+1$.
Using
inequalities
(\ref{21.10.2009_1}), (\ref{21.10.2009_2}),
we obtain
$$
                    \sigma_f (B_1,\dots,B_l; C)
                        =
                            \sum_{\upsilon\in S} |\sigma_{f,\upsilon} (B_1,\dots,B_l; C)|
                                \le
                                    \sigma'_1 + \sigma''_1
                                        +
                                            \sigma_f (B_1,\dots,B_l; \t{C})
                                                \le
$$
\begin{equation}\label{21.10.2009_3}
    \le
        2^L 3^{lm} \eps^{1/4} N^{l-m} + \eta \tau_f (B'_1,\dots,B'_l ; C)
            +
                \sigma_f (B_1,\dots,B_l; \t{C}) \,.
\end{equation}
Recalling that
$\eta = |S|^{-1} \eps_1/16$, using
condition (\ref{cond:eps_small})
and formulas (\ref{tmp:10.12.2009_1}), (\ref{21.10.2009_3}),
we get
\begin{equation}\label{tmp:25.10.2009_1}
    \sigma_{f} (B_1,\dots,B_l; C)
            \le
                \sigma_{f} (B_1,\dots,B_l; \t{C})
                            +
                                2^{-2} \eps_1 |S|^{-1} \tau_f (B'_1,\dots,B'_l ; C) \,.
\end{equation}
Finally,
$$
    \tau_f (B_1,\dots,B_l; \t{C}) + \tau_f (B'_1,\dots,B'_l; C)
        =
            \tau_f (B_1,\dots,B_l; C)
$$
and by (\ref{f:number_of_good_tuples_below}), we obtain
(\ref{f:solutions_DONOT_decrease}).
This concludes the proof of the lemma.

Now return to the proof of the theorem.
We use Lemma \ref{l:density_inc} inductively.
At zeroth step, we have
$\a \ge |S|^{-1} \eps_1$
and $\o$ is any positive number such that $\o \le \eps_1 |S|^{-1}$.
Suppose that our algorithm was applied $h$ times.
Thus we obtain the sets $C_j$, $j\in [h]$,
every $C_j$ is a set of level $j$ and two
sequences of $\zeta_1,\dots,\zeta_h$,  $\o_1,\dots,\o_h$,
each $\zeta_j$ depends on $\o_j$.
Using inequality (\ref{f:sigma_a+}),
we see that for any $j$ the following holds
$$
    \sigma_f (B_1,\dots,B_l ; \Gr^l)
        \le
            \sigma_f (B_1,\dots,B_l ; C_j)
                +
                    2^{-2} \eps_1 |S|^{-1} \tau_f (B_1,\dots,B_l ; \Gr^l)
                        \le
$$
$$
                        \le
    \sigma_f (B_1,\dots,B_l ; C_j)
       +
    2^{-2} \eps_1 N^{l-m} \,.
$$
So, we can assume that for all $j\in [h]$,
we have $\sigma_f (B_1,\dots,B_l ; C_j) \ge \eps_1 (2|S|)^{-1} \tau_f (B_1,\dots,B_l ; C_j)$.
Further, in view of (\ref{f:sigma_le_tau}),
we can suppose that for all $j\in [h]$
the following holds $\tau_f (B_1,\dots,B_l ; C_j) \ge 2^{-1} \eps_1 N^{l-m}$.
Thus
for any $j\in [h]$, we can take
$\omega_j := 2^{-1} \eps_1 |S|^{-1}$
and hence for all $j\in [h]$ the following holds
$\zeta_j \ge \l( \frac{\eps^2_1}{64(1+2 \eps^{-1}_1 |S|)|S|^2} \r)^{lm  2^{3ml}} = \zeta_*$.
By (\ref{f:solutions_DONOT_decrease}), we get
$$
    \tau_f (B_1,\dots,B_l ; C_j) \le (1- \zeta_*)^j \tau_f (B_1,\dots,B_l ; \Gr^l)
        \le
            (1- \zeta_*)^j |S| N^{l-m}
$$
and
our algorithm must stop after at most
$2^4 \log (2|S| \eps^{-1}_1) \l( \frac{\eps^2_1}{64(1+2\eps^{-1}_1 |S|)|S|^2} \r)^{-lm  2^{3ml}} := h$
number of steps.
Clearly, the maximal level of any set which appears in our induction procedure
does not exceed $h$.
Using the inequality $\tau_{f} (B_1,\dots,B_l ; C_h) \ge 2^{-1} \eps_1 N^{l-m}$,
we get from (\ref{cond:eps_small}) the dependence between parameters $\eps$ and $\eps_1$
\begin{equation}\label{f:last1}
    \eps
        \le
            2^{-4h-20} |S|^{-8}
                \eps_1^{8} \l( \frac{\eps_1^2}{64(1+2\eps_1^{-1}|S|)|S|^2} \r)^{4lm 2^{3ml}}
\end{equation}
Besides, by (\ref{cond:N_new}), we obtain
\begin{equation}\label{f:last2}
    N \ge 2^{2l+1} \l( \frac{64(1+2\eps_1^{-1} |S|)|S|^2}{\eps_1^2} \r)^{lm 2^{3ml}} |S| \eps^{-1}_1\,.
\end{equation}
If (\ref{f:last1}) and (\ref{f:last2}) hold then we have a contradiction.
Small calculation shows that the last two inequalities imply (\ref{concl:sigma_small}).
This completes the proof of Theorem \ref{t:main_skew_norms}.

\Note
Probably, the result above suggests that
we have
a phenomenon
in
spirit of
the
well--known sum--product phenomenon (see e.g. \cite{ESz} or \cite{BKT})
or the dichotomy phenomenon (see e.g. \cite{Green_finite_fields}).
Indeed, Theorem \ref{t:main_skew_norms} is trivial in two opposite cases :
if $a(x)$ is a linear function then $\sigma_f (B_1,\dots, B_l; \Gr)$
is small by (\ref{cond:rectangle_small})
and
if $a(x)$ is far from all linear functions (e.g. $a(x)=x^2$) then (\ref{concl:sigma_small})
is small by (\ref{f:sigma_le_tau}) and
an appropriate upper bound for $\tau_f (B_1,\dots, B_l; \Gr)$.
So,
in some sense
our result
can be interpreted
that there is no function
linear and non--linear
simultaneously.

\refstepcounter{section}
\label{sec:applications}

\begin{center}
{\bf \arabic{section}. Applications.}
\end{center}

First of all let us note a simple property of Gowers norms
(see e.g. \cite{Gow_m}).

\Lemma
{\it
    Let $d\ge 2$ be an integer, and $f:\Gr^k \to \C$ be a function.
    Let also $u_1, \dots, u_d : \Gr^d \to [-1,1]$ be any functions such that
    $u_i (\v{x}) = u_i ((\v{x})_{(i)})$, $i=1,\dots,d$, $\v{x} = (x_1,\dots,x_d)$.
    Then
    $$
        \| f(\v{x}) e( \prod_{i=1}^d u_i (\v{x}) ) \|_{U^d} = \| f(\v{x}) \|_{U^d} \,.
    $$
}
\label{l:exp_in_G_norms}

Also, we need in a convexity lemma.

\Lemma
{\it
    Let $\kappa>0$ and $h(x) = 1/(\log x)^\kappa$.
    Then for any real numbers $x_1,\dots,x_N > 1$, we have
    $$
        \frac{1}{N} \sum_{i=1}^N h(x_i) \le h \l( \frac{1}{N} \sum_{i=1}^N x_i \r) \,.
    $$
}
\label{l:convexity_1/log}

Now
we
obtain a corollary of Theorem \ref{t:main_skew_norms}.
Clearly, the result below implies Theorem \ref{t:intr_rect}.

\Cor
{\it
    Let $\eps\in (0,1)$ be a real number,
    $a:\Gr \to \Gr$ be a function, and
    $f (x,y) = e(x \cdot a(y))$.
    Suppose that
    \begin{equation}\label{cond:cor_f_rect}
        \| f(x,y) \|_{U^2} \le \eps \,.
    \end{equation}
    Then
    \begin{equation}\label{f:eps_t^2}
        \| f(x, x+y) \|_{\U^{2}}
            \le
                3^4 \max\l\{ 2^{30} \l( \frac{2^8}{\log (1/\eps)} \r)^{2^{-21}},
                        4 \l( \frac{2^{8}}{N} \r)^{2^{-18}} \r\}
                 N^{4}
        \,.
    \end{equation}
    If $t$ is a positive integer and
    \begin{equation}\label{cond:cor_f_rect+}
        \| f(x,y_1+\dots+y_t) \|_{U^{t+1}} \le \eps
    \end{equation}
    then
    \begin{equation}\label{f:eps_t^2+}
        \| f(x, x+y_1+\dots+y_t) \|_{\U^{t+1}}
            \ll_{t} \max \l\{ \frac{1}{\log ( 1/\eps)}, N^{-c(t)} \r\} N^{2t+2}
        \,,
    \end{equation}
    where $c(t)>0$ some constant depends on $t$ only.
}
\label{cor:rectangle_imply}

\Note
{
Clearly,
formula (\ref{f:eps_t^2}) of the corollary above
holds
for
any
function $f_\la (x,y) = e((x-\la) \cdot a(y))$,
where $\la \in \Gr$ is an arbitrary element.
}
\label{note:shift_la}

\Proof
For $w=(w_1,\dots,w_d)$, $d\in \N$ and $\o \in \{0,1\}^d$
we write $w_\o$ for $\sum_{i=1}^d w^\o_i$.
Any sequence of points $(w_\o)_{\o \in \{0,1\}^d} \in \Gr$
is called {\it $d$--dimensional cube} (see \cite{Gow_m}).
So, we numerate the points of an arbitrary  $d$--dimensional cube by index $\o \in \{0,1\}^d$.

First of all let us note that
the dimension of all $d$---dimensional cubes equals $d+1$
and
any $d$---dimensional
cube
$(z_1,\dots,z_{2^d})$
is a solution of a system of full rank
having
$\sum_{i=0}^{d-2} \binom{d}{i} = 2^d-d-1$
linear equations, say
\begin{displaymath}
\left\{ \begin{array}{ll}
        \sum_{\omega \in \{0,1\}^d} z_{\omega} \cdot (-1)^{|\omega|} = 0 \,, \\
        \sum_{\omega \in \{0,1\}^{d}} z_{\omega} \cdot (-1)^{|\omega|} = 0 \,,
            \quad \o_i = 0,\, i\in [d] \,. \\
        \dots \\
        \sum_{\omega \in \{0,1\}^{d}} z_{\omega} \cdot (-1)^{|\omega|} = 0 \,,
            \quad \o_{i_1} = \dots = \o_{i_{d-2}} = 0,\, i_j \in [d] \,. \\
\end{array} \right.
\end{displaymath}
By $S$ denote the last system.
The fact that
$S$
has rank equals $2^d-d-1$ can
be obtained by successfully considering of its equations
from the first to the last (see variables $z_{1,\dots,1}$,
$z_{0,1,\dots,1}, \dots, z_{1,\dots,1,0}$ and so on).

Our task is count the quantity
$$
    \sigma_t := \| f(x, x+y_1+\dots+y_t) \|_{\U^{t+1}} =
        \sum_{x,x'}\, \sum_{y_1,\dots,y_t}\, \sum_{y'_1,\dots,y'_t}
            e( \sum_{\o \in \{0,1\}^{t+1}} x_{\o_1} \cdot a(x_{\o_1} + y_{(\o_2,\dots,\o_{t+1})}) ) \,,
$$
where $y=(y_1,\dots,y_t)$.
Let us change the variables $x \to x-y_1$, $x' \to x'-y_1$, $y_1'-y_1 = \Delta$.
Let $z=(x,\Delta,y_2,\dots,y_t)$,
$\o \in \{0,1\}^{t+1}$, $\o = (\bar{\o}_1,\eta,\bar{\o}_2)$,
where $\bar{\o}_1,\eta \in \{ 0,1 \}$, $\bar{\o}_2 \in \{0,1\}^{t-1}$.
Put $z_\o = x_{\bar{\o}_1} + \eta \Delta + y_{\bar{\o}_2}$.
So, $z_\o$ depends on $x,\Delta,y$.
We write $\eta=\eta(z_\o)$.
Clearly,
$$
    \sigma_t
        =
            \sum_{x,x',\Delta}\,\, \sum_{y_2,\dots,y_t,\,y'_2,\dots,y'_t}\,
                e ( \sum_{\o \in \{0,1\}^{t+1}} x_{\bar{\o}_1} \cdot a(x_{\bar{\o}_1}
                    + \eta \Delta + y_{(\o_2,\dots,\o_t)} )
    \m
$$
\begin{equation}\label{tmp:23_10_2009_1}
    \m
                    N \d_0 (\sum_{\o \in \{0,1\}^{t+1}} (-1)^{|\o|} a(z_\o)) \,,
\end{equation}
where $\d_0 (x) = 1$ if $x=0$ and zero otherwise.
It is easy to see that the sequence $(z_\o)_{\o \in \{0,1\}^{t+1}}$
form a cube.
So, our sum $\sigma_t$ has form (\ref{f:sigma_form}).
To prove Corollary \ref{cor:rectangle_imply} we need to check
condition (\ref{cond:rectangle_small}) but firstly
let us prove the corollary in the simplest case $t=1$.
In the situation formula (\ref{tmp:23_10_2009_1}) is
$$
    \sigma_1
        =
            \sum_{x,x',\Delta}
                e( xa(x) - x' a(x') - x a(x+\Delta) + x' a(x'+\Delta) )
                    \cdot N \d_0 ( a(x) - a(x') - a(x+\Delta) + a(x'+\Delta) )
$$
$$
    =
        \sum_{x,x',\Delta}
            e ( (x-x')(a(x) - a(x+\Delta)) )
                \cdot N \d_0 ( a(x) - a(x') - a(x+\Delta) + a(x'+\Delta) )
$$
and condition (\ref{cond:rectangle_small}) holds because (\ref{cond:cor_f_rect})
and
$S$ contains just two equations
in the case.
Thus, we prove (\ref{f:eps_t^2}), provided by $t=1$.

Now
suppose that $t\ge 2$.
Let $\sum_\o a(z_\o) \bar{\eta}_\o = 0$ be an equation from $S$.
Suppose that there is $\bar{\eta}_\o \neq 0$ such that $\eta(z_\o) = 1$
(see definition of $z_\o$).
Without losing of generality, one can suppose that $\bar{\eta}_{\v{1}} = -1$,
where $\v{1} = (1,\dots,1)$.
Then
$a(z_{\v{1}}) = \sum_{\o,\, \o \neq \v{1}} a(z_\o) \bar{\eta}_\o$
and we can substitute it into formula (\ref{tmp:23_10_2009_1}).
After that change the variables $x\to x-y_2$, $x'\to x'-y_2$,
$y'_2-y_2 = q$.
Consider the variables $x',q,y'_3,\dots,y'_t$.
We see there are terms
$e(x' a(x+q+y'_3+\dots+y'_t))$,
$e(x' a(x+\Delta+q+y'_3+\dots+y'_t))$
and
may be
the term $e((x'-y_2) a(x'+q+y'_3+\dots+y'_t))$,
containing $x',q,y'_3,\dots,y'_t$
but there is no term
$e((x'-y_2) a(x'+\Delta+q+y'_3+\dots+y'_t))$
in the expression for $\sigma_t$.
Put
\begin{equation}\label{f:Theta}
    \Theta^{(1)}_{y_2,x,\Delta} (x',q,y'_3,\dots,y'_t)
        =
    e( x' a(x+q+y'_3+\dots+y'_t) + x' a(x+\Delta+q+y'_3+\dots+y'_t))
\end{equation}
and
\begin{equation}\label{f:Theta_2}
    \Theta^{(2)}_{y_2,x,\Delta} (x',q,y'_3,\dots,y'_t)
        =
            \Theta^{(1)}_{y_2,x,\Delta} (x',q,y'_3,\dots,y'_t)
                \cdot
    e( (x'-y_2) a(x'+q+y'_3+\dots+y'_t))
\end{equation}
provided by the term $e( (x'-y_2) a(x'+q+y'_3+\dots+y'_t))$ exists.
%
%
%
Clearly, another multiples do not contain
these variables.
Using Lemma \ref{l:exp_in_G_norms} and summation over $y_2$, we get
$$
    \sum_{y_2,x,\Delta} \| \Theta^{(2)}_{y_2,x} (x',q,y'_3,\dots,y'_t) \|_{\U^t}
        \le
            N^2 \| f(x,y_1+\dots+y_t) \|_{\U^{t+1}} \,.
$$
Similarly, summing over $x,\Delta$, we have
$$
   \sum_{y_2,x,\Delta} \| \Theta^{(1)}_{y_2,x} (x',q,y'_3,\dots,y'_t) \|_{\U^t}
        \le
            N^2 \| f(x,y_1+\dots+y_t) \|_{\U^{t+1}} \,.
$$
By Lemma \ref{l:convexity_1/log}, we obtain (\ref{f:eps_t^2+}).
Now suppose that
for our equation $\sum_\o a(z_\o) \bar{\eta}_\o = 0$ from $S$
for any
$\bar{\eta}_\o \neq 0$ we have $\eta(z_\o) = 0$.
But changing the variables
$y_2 \to y_2 - \Delta$, $y'_2 \to y'_2 - \Delta$
we reducing our expression to the previous case.
This concludes the proof.

\Note If we take  $a(x)= {\rm const}$ then, clearly, (\ref{f:eps_t^2}) does not hold.
On the other hand it is easy to see that (\ref{cond:cor_f_rect})
is also has no place.
If we put $a(x)=x$  then (\ref{f:eps_t^2+}) and (\ref{cond:cor_f_rect+})
do not hold for $t=2$.
For $t\ge 2$ take $a(x)=x^{t-1}$.
Then the norm of the functions from (\ref{cond:cor_f_rect+})
and (\ref{f:eps_t^2+})
equal one.
Thus, in some sense our conditions (\ref{cond:cor_f_rect}), (\ref{cond:cor_f_rect+})
are necessary.

\Note
Let us analyze a little bit more
general functions of the form
$f(z_1+\dots+z_r+x_1+\dots+x_s, x_1+\dots+x_s+y_1+\dots+y_t)$, $r$ is a positive integer,
$s,t\ge 0$.
If $r\ge 2$ then $\| f \|_{U^{s+t+r}} = 1$ for any function $a(x)$.
Further, if $r=0$, $t=0$ and $\Gr = \Z_p$, $p$ is a  prime number then put $a(x) = x^{-1}$.
It is easy to see that inequality (\ref{cond:cor_f_rect}) takes place
(e.g. see Corollary \ref{cor:injection})
and $\| f \|_{U^{s}} = 1$.
The case $r=s=0$ is not interesting and consider the situation
$r=0$, $s\ge 2$, $t\ge 1$.
It is easy to see that $\| e(x a(y_1+\dots+y_{s+t})\|_{U^{s+t+1}} = o(1)$
trivially implies
$\| f(x_1+\dots+x_s,x_1+\dots+x_s+y_1+\dots+y_t) \|_{U^{s+t}} = o(1)$
(see the first step of proof of Corollary \ref{cor:rectangle_imply})
and for $a(x)=x^{s+t-2}$,
we have $\| f(x_1+\dots+x_s,x_1+\dots+x_s+y_1+\dots+y_t) \|_{U^{s+t}} = 1$.
Finally, in the case
$r=1$, $s,t\ge 1$ we can apply the same argument.
Thus,
the
choice
$r=0$, $s=1$, $t\ge 1$ is the only situation when
we
can obtain non--trivial upper bounds (under some analogs of assumption (\ref{cond:cor_f_rect+}))
for $f(z_1+\dots+z_r+x_1+\dots+x_s, x_1+\dots+x_s+y_1+\dots+y_t)$.


\Cor
{\it
    Let $K$ be a positive integer,
    $a:\Gr \to \Gr$ be a function, and
    $f (x,y) = e(x \cdot a(y))$.
    Suppose that
    \begin{equation}\label{cond:preimage}
        | \{ y\in \Gr ~:~ a(y) = j \}| \le K\,, \quad j\in \Gr \,.
    \end{equation}
    Then
    $$
    \| f(x, x+y) \|_{\U^{2}}
            \le
                3^4 \max\l\{ 2^{30} \l( \frac{2^8}{\log (1/\eps_*)} \r)^{2^{-19}},
                        4 \l( \frac{2^{8}}{N} \r)^{2^{-16}} \r\}
                 N^{4}
                \,,
    $$
    where $\eps_* = (K/N)^{1/4}$.
}
\label{cor:injection}
\\
\Proof
It is easy to see that
$$
    \| f(x,y) \|^4_{U^2} = \frac{1}{N^4} \sum_{x,x',y,y'} e( (x-x') (a(y)-a(y')))
        =
            \frac{1}{N^2} \sum_{j\in \Gr} |M_j|^2 \,,
$$
where $M_j = \{ y\in \Gr ~:~ a(y) = j \}$.
By assumption, we have $|M_j| \le K$ for all $j\in \Gr$.
Besides $\sum_j |M_j| = N$.
Hence $\|f(x,y) \|_{U^2} \le (K/N)^{1/4}$ and the statement follows from
Corollary \ref{cor:rectangle_imply}.
This completes the proof.

Thus, if, say, $K=N^{1-\kappa}$, $N\to \infty$, $\kappa>0$
then
$\| f(x, x+y) \|_{U^{2}} \ll 1/(\log N)^c$,
where $c=c(\kappa)$.
Even in the situation when preimage  of any (or almost any) point $j\in \Gr$
has the cardinality $o(N)$ we have
a non--trivial bound for $\| f(x, x+y) \|_{U^{2}}$.

\Note
It is easy to see that
$$
    \| e(x a(y_1+\dots+y_t)) \|_{\U^{t+1}}
        =
            N^2 \cdot \# \{ y_1,\dots,y_t,y'_1,\dots,y'_t
                        ~:~ \sum_{\o \in \{0,1\}^t} a(y^\o) (-1)^{|\o|} = 0 \} \,.
$$
So, condition (\ref{cond:cor_f_rect+}) can be interpreted
as the requirement for $a(x)$ to be far from "polynomial of degree $t-1$"
(or $(t-1)$--step nilsequence, more precisely, see e.g. \cite{GT_U3,H-Kra,Ziegler_1}).

Finally, we apply Corollary \ref{cor:injection} to
a family of subsets of $\Gr \m \Gr$, $\Gr = \Z_N$.
Let us recall some formulas from Fourier analysis.
Let $f:\Gr \to \C$ be a function.
By $\F{f}(\xi)$ denote the Fourier transformation of $f$
\begin{equation}\label{F:Fourier}
  \F{f}(\xi) =  \sum_{x \in \Gr} f(x) e( -\xi \cdot x) \,,
\end{equation}
We have
\begin{equation}\label{F_Par}
    \sum_{x\in \Gr} |f(x)|^2
        =
            \frac{1}{N} \sum_{\xi \in \F{\Gr}} |\widehat{f} (\xi)|^2
                \quad \quad \mbox{ (Parseval identity) }
\end{equation}
and
\begin{equation}\label{f:inverse}
    f(x) = \frac{1}{N} \sum_{\xi \in \F{\Gr}} \F{f}(\xi) e(\xi \cdot x)
        \quad \quad \mbox{ (the inverse formula) }\,.
\end{equation}
If
$$
    (f*g) (x) := \sum_{y\in \Gr} f(y) g(x-y)
$$
 then
\begin{equation}\label{f:F_svertka}
    \F{f*g} = \F{f} \F{g} \,.
\end{equation}

Now we can formulate our corollary.

\Cor
{\it
    Let $\Gr = \Z_N$, $N$ be a prime number.
    Let also $a:\Gr \to \Gr$ be a function,
    $a(y) \neq 0$ for all $y\in \Gr$.
    Suppose that $P$ is an arithmetic progression, $|P| \ge 2^5 N^{3/4}$,
    $f(x,y) = P(x \cdot a(y)) - |P| / N$,
    and
    inequality (\ref{cond:preimage}) holds with $K\le 2^{-18} |P|^4 / N^3$.
    Then
    \begin{equation}\label{f:cor_sets_P_bound}
     \| f(x, x+y) \|_{\U^{2}}
            \le
                3^4 \max\l\{ 2^{30} \l( \frac{2^8}{\log (1/\eps_*)} \r)^{2^{-21}},
                        4 \l( \frac{2^{8}}{N} \r)^{2^{-18}} \r\}
                 N^{4}
                \,,
    \end{equation}
    where $\eps_* = (4K/N)^{1/48}$.
}
\label{cor:sets_P}
\\
\Proof
To reduce some logarithms in our bounds we use a well--known trick.
Without losing of generality, suppose that
$P = \{ 0,1,\dots,|P|-1 \}$.
Let $P_1 = \{ 0,1,\dots, t-1 \}$, $t= [ |P|^{2/3} K^{1/12} N^{1/4} 2]$,
and $W(x) = |P_1|^{-1} (P*P_1)(x)$.
Clearly, $0\le W(x) \le 1$, $\sum_x W(x) = |P|$ and $W(x) = P(x)$
for all but at most $2t$ points $x\in \Gr$.
By Parseval, formula (\ref{f:F_svertka}) and the Cauchy--Schwartz, we obtain
\begin{equation}\label{tmp:25.10.2009_3}
    \sum_{r} |\F{W} (r)|
        \le
            \frac{1}{|P_1|} \sqrt{N |P_1|} \sqrt{N |P|}
                =
                    N |P|^{1/2} t^{-1/2} \,.
\end{equation}
By the triangle inequality for $\U^2$---norm and (\ref{f:inverse}), we get
$$
    \| f \|_{\U^2} \le \| W - |P|/N \|_{\U^2} + 32 \cdot t N^{3}
        =
$$
$$
        =
            \frac{1}{N^4} \sum_{(r_1,\dots,r_4) \neq \v{0}} \F{W} (r_1) \dots \F{W} (r_4)\,
                \sum_{x,x',y,y'} e(r_1 x a(y)) e(-r_2 x' a(y)) e(-r_3 x a(y')) e(r_4 x' a(y'))
                    +
$$
\begin{equation}\label{tmp:25.10.2009_4}
                    +
                        32 \cdot t N^{3} \,.
\end{equation}
Clearly, $\| e(r\cdot x a(y)) \|_{U^2} = \| e(x a(y)) \|_{U^2}$, provided by $r\neq 0$.
We have $\| e(x a(y)) \|_{U^2} \le (K/N)^{1/4}$
(see the proof of Corollary \ref{cor:injection}).
Using the last inequality, Lemma \ref{l:U_k_sense}, formulas (\ref{tmp:25.10.2009_3}),
(\ref{tmp:25.10.2009_4}), we obtain
$$
    \| f \|_{\U^2}
        \le
            \l( \frac{|P|}{t} \r)^{2}
                \l( \frac{K}{N} \r)^{1/4} N^3
                    +
                        32 \cdot t N^3
                            \le
                                2^{8} \l( \frac{K}{N} \r)^{1/12} N^4 \,.
$$
Using Corollary \ref{cor:rectangle_imply}, we get (\ref{f:cor_sets_P_bound}).
This concludes the proof.

\Note There is another alternative way to obtain
the smallness of $\| \cdot \|_{U^2}$---norm
of the function $f(x,y) = P(x\cdot a(y)) - |P|/N = P(x\cdot y) - |P|/N$
in the situation $\Gr = \Z_N$
and
$a(y)=y$.
It is easy to see,
using multiplicative characters,
that the $\| f \|_{U^2}$ is small by
P\'{o}lya---Vinogradov  inequality
(see e.g. \cite{Vinogradov_book}).


\end{document}